\documentstyle[12pt]{article}
\newtheorem{lemma}{Lemma}[section]
\newtheorem{theorem}[lemma]{Theorem}

\newtheorem{claim}[lemma]{Claim}

\newtheorem{proposition}[lemma]{Proposition}

\newtheorem{conjecture}[lemma]{Conjecture}

\newcommand{\comment}[1]{}

\newcommand{\text}[1]{\quad\mbox{#1}\quad}
\def\beq{\begin{equation}}\def\eeq{\end{equation}}
\def\beqn{\begin{eqnarray}}\def\eeqn{\end{eqnarray}}
\def\pont{\hspace{-6pt}{\bf.\ }}
\def\eps{\varepsilon}
\def\qed{\ifhmode\unskip\nobreak\fi\quad\ifmmode\Box\else$\Box$\fi}

\title{Partitioning $2$-edge-colored graphs by monochromatic paths and
cycles}

\author{J\'ozsef Balogh\thanks{Research  supported in part   by NSF CAREER
Grant DMS-0745185, UIUC Campus Research Board Grant 11067, and OTKA
Grant K~76099.}\\
\small Department of Mathematical Sciences\\[-0.8ex]
\small University of Illinois\\[-0.8ex]
\small Urbana, IL 61801\\[-0.8ex]
\small \texttt{jobal@math.uiuc.edu}\\
J\'anos Bar\'at
\thanks{Research is supported by OTKA Grants PD~75837 and K~76099,
and the J\'anos Bolyai Research Scholarship of the Hungarian Academy
of Sciences. Present address: School of Mathematical
Sciences, Monash University, 3800 Victoria.}\\
\small Department of Computer Science and Systems Technology\\[-0.8ex]
\small University of Pannonia, Egyetem u.\ 10, 8200 Veszpr\'em,
Hungary\\[-0.8ex]
\small \texttt{barat@dcs.vein.hu}\\
D\'aniel Gerbner\\
\small Hungarian Academy of Sciences, Alfr\'ed R\'enyi Institute\\[-0.8ex]
\small of Mathematics, P.O.B. 127, Budapest H-1364, Hungary\\[-0.8ex]
\small \texttt{gerbner@renyi.hu}\\
Andr\'as Gy\'arf\'as\thanks{Research supported in part by OTKA Grant K104373.}\\
\small Alfr\'ed R\'enyi Institute of Mathematics\\[-0.8ex]
\small Hungarian Academy of Sciences\\[-0.8ex]
\small Budapest, P.O. Box 127\\[-0.8ex]
\small Budapest, Hungary, H-1364 \\[-0.8ex]
\small \texttt{gyarfas.andras@renyi.mta.hu}\\
and\\
G\'{a}bor N. S\'ark\"ozy\thanks{Research supported in part by
NSF Grant DMS-0968699 and by OTKA Grant K104373.}\\
\small Computer Science Department\\[-0.8ex]
\small Worcester Polytechnic Institute\\[-0.8ex]
\small Worcester, MA, USA 01609\\[-0.8ex]
\small \texttt{gsarkozy@cs.wpi.edu}\\[-0.8ex] }
\begin{document}
\maketitle

\begin{abstract}
We present results on partitioning the vertices of  $2$-edge-colored
graphs into monochromatic paths and cycles. We prove asymptotically
the two-color case of a conjecture of S\'ark\"ozy: the vertex set of
every $2$-edge-colored graph can be partitioned into at most
$2\alpha(G)$ monochromatic cycles, where $\alpha(G)$ denotes the
independence number of $G$. Another direction, emerged recently from
a conjecture of Schelp, is to consider colorings of graphs with
given minimum degree. We prove that apart from $o(|V(G)|)$ vertices,
the vertex set of any $2$-edge-colored graph $G$ with minimum degree
at least ${(1+\eps)3|V(G)|\over 4}$ can be covered by the vertices
of two vertex disjoint monochromatic cycles of distinct colors.
Finally, under the assumption that $\overline{G}$ does not contain a
fixed bipartite graph $H$, we show that in every $2$-edge-coloring
of $G$, $|V(G)|-c(H)$ vertices can be covered by two vertex disjoint
paths of different colors, where $c(H)$ is a constant depending only
on $H$. In particular, we prove that $c(C_4)=1$, which is best
possible.\footnote{Part of the research reported in this paper was
done at the 3rd Eml\'ekt\'abla Workshop (2011) in Balatonalm\'adi,
Hungary.}
\end{abstract}

\section{Background, summary of results.}

In this paper, we consider some conjectures about partitioning
vertices of edge-colored graphs into monochromatic cycles or paths.
For simplicity, colored graphs means edge-colored graphs in this
paper. In this context it is conventional to accept {\it empty
graphs and one-vertex graphs} as a path or a cycle (of any color)
and also {\it any edge} as a path or a cycle (in its color). With
this convention one can define the {\em cycle (or path) partition
number} of any colored graph $G$ as the minimum number of vertex
disjoint monochromatic cycles (or paths) needed to cover the vertex
set of $G$. For complete graphs, \cite{EGP} posed the following
conjecture.

\begin{conjecture}\pont\label{egypconj} The cycle partition number of
any $t$-colored complete graph $K_n$ is $t$.
\end{conjecture}

The $t=2$ case of this conjecture was stated earlier by Lehel in a
stronger form, requiring that the colors of the two cycles must be
different.  After some initial results~\cite{AY, GY},
 \L uczak,  R\"{o}dl and Szemer\'edi~\cite{LRS}
 proved Lehel's conjecture for large enough $n$, which can be
considered as a birth of certain advanced applications of the
Regularity Lemma. A more elementary proof, still for large enough
$n$, was obtained by Allen~\cite{A}. Finally, Bessy and Thomass\'e
\cite{BT} found a completely elementary inductive proof for every
$n$.

The $t=3$ case of Conjecture \ref{egypconj} was solved
asymptotically in~\cite{GRSSZ}. Pokrovskiy~\cite{PR} showed recently
(with a nice elementary proof) that the path partition number of any
$3$-colored $K_n$ is at most three (for any $n\ge 1$). But then
surprisingly Pokrovskiy \cite{PR1} found a counterexample to
Conjecture \ref{egypconj} for all $t\geq 3$. However, in the
counterexample all but one vertex can be covered by $t$ vertex
disjoint monochromatic cycles.

For general $t$, the best bound for the cycle partition number is
$O(t\log{t})$, see \cite{rlogr}. Note that it is far from obvious
that the cycle partition number of $K_n$ can be bounded by {\it any}
function of $t$.

We address the extension of the cycle and path partition numbers
from complete graphs to arbitrary graphs $G$. If we want these
numbers to be independent of $|V(G)|$, some other parameter of $G$
must be included. We consider three of these parameters.

Let $\alpha(G)$ denote the independence number of $G$, the maximum
number of pairwise non-adjacent vertices of $G$. The role of
$\alpha(G)$ in results on colorings of non-complete graphs was
observed in \cite{GST,Gydisz,GYS1} and in S\'ark\"ozy \cite{SA} who
extended Conjecture ~\ref{egypconj} to the following.

\begin{conjecture}\pont\label{saconj}
The cycle partition number of any $t$-colored graph $G$ is
$t\alpha(G)$.
\end{conjecture}

For $t=1$, Conjecture~\ref{saconj} is a well-known result of P\'osa
\cite{posa} (and clearly best possible). For $t=2$ it is also best
possible, shown by vertex disjoint copies of triangles, each colored
using two colors. To prove Conjecture~\ref{saconj} for $t=2$ and
arbitrary $\alpha(G)$ seems very difficult (considering the
complexity of the proof for $\alpha(G)=1$ in \cite{BT}). Then again
the counterexample of Pokrovskiy \cite{PR1} shows that the
conjecture is not true in this form for any $t\geq 3$. Perhaps the
following weakening of the conjecture is true.

\begin{conjecture}\pont\label{saconj1}
Let $G$ be a $t$-colored graph with $\alpha(G)=\alpha$. Then there
exists a constant $c=c(\alpha,t)$ such that $t\alpha$ vertex
disjoint monochromatic cycles of $G$ cover at least $n-c$ vertices.
\end{conjecture}

Pokrovskiy's example implies that $c\geq \alpha$ must be true. We
cannot prove this conjecture even for $t=2$, we can only prove the
following weaker asymptotic result.

\begin{theorem}\pont\label{saass}
For every positive $\eta$ and $\alpha$, there exists an
$n_0(\eta,\alpha)$ such that the following holds. If $G$ is a
2-colored graph on $n$ vertices, $n\geq n_0$, $\alpha(G)=\alpha$,
then there are  at most $2\alpha$ vertex disjoint monochromatic
cycles covering at least $(1-\eta)n$ vertices of $V(G)$.
\end{theorem}

Recently, Schelp \cite{SCH} suggested in a posthumous paper to
strengthen certain Ramsey problems from complete graphs to graphs of
given minimum degree. In particular, he conjectured that with
$m=R(P_n,P_n)$, minimum degree ${3m\over 4}$ is sufficient to find a
monochromatic path $P_n$ in any $2$-colored graph of order
$m$.\footnote{Some progress towards this conjecture have been done
in \cite{GYS2} and \cite{BBLS}.} Influenced by this conjecture, here
we pose the following conjecture.

\begin{conjecture}\pont\label{schconj}
If $G$ is an $n$-vertex graph with $\delta(G)>3n/4$
then in any $2$-edge-coloring of $G$, there are two vertex disjoint
monochromatic cycles of different colors, which together cover
$V(G)$.
\end{conjecture}

That is, the above mentioned Bessy-Thomass\'e result \cite{BT} would
hold for graphs with minimum degree larger than $3n/4$. Note that
the condition $\delta(G)\ge {3|V(G)|\over 4}$ is sharp. Indeed,
consider the following $n$-vertex graph, where $n=4m$. We partition
the vertex set into four parts $A_1,A_2,A_3,A_4$ with $|A_i|=m$.
There are no edges from $A_1$ to $A_2$ and from $A_3$ to $A_4$.
Edges in $[A_1,A_3],[A_2,A_4]$ are red and edges in
$[A_1,A_4],[A_2,A_3]$ are blue, inside the classes any coloring is
allowed. In such an edge-colored graph, there are no two vertex
disjoint monochromatic cycles of {\em different colors} covering
$G$, while the minimum degree is $3m-1={3n\over 4}-1$.

We prove Conjecture~\ref{schconj} in the following asymptotic sense.

\begin{theorem}\pont\label{mindegass}
For every $\eta > 0$, there is an $n_0(\eta)$ such that the
following holds. If $G$ is an $n$-vertex graph with $n\geq n_0$ and
$\delta(G)>(\frac{3}{4}+\eta) n$, then every $2$-edge-coloring of
$G$ admits two vertex disjoint monochromatic cycles of different
colors covering at least $(1-\eta)n$ vertices of $G$.
\end{theorem}

The proofs of Theorems~\ref{saass} and \ref{mindegass} follow a
method of \L uczak~\cite{L}.
The crucial idea is that the words ``cycles" or ``paths" in a statement to be proved
are replaced by the words ``connected matchings". In a {\it connected matching},
the edges of the matching are in the same component of the graph.\footnote{When the
edges are colored, a connected red matching is a matching in a red
component.} We prove first this weaker result, then we   apply to the cluster graph of a
regular partition of the target graph. Through several technical
details, the regularity of the partition is used to ``lift back" the
connected matching of the cluster graph to a path or cycle in the
original graph. In our case, the relaxed versions of
Theorems~\ref{saass} and \ref{mindegass} for connected matchings are
stated and proved in Section \ref{matchingresults} (Theorem
\ref{connmatch} and \ref{degmatch}).

Another possibility to extend Conjecture \ref{egypconj} to more
general graphs is to consider a graph $G$, whose complement does not
contain a fixed bipartite graph $H$. This brings in a different
flavor, since these graphs are very dense, they have ${|V(G)|\choose
2}-o(|V(G)|^2)$ edges. In return, we  prove sharper results in this
case. We also state a more general conjecture.

\begin{conjecture}\pont\label{conj2} Let $H$ be a graph with chromatic
number $k+1$ and let $G$ be an $t$-edge-colored graph on $n$
vertices such that $H$ is not a subgraph of $\overline{G}$. Then
there exists a constant $c=c(H,k,t)$ such that $kt$ vertex disjoint
monochromatic paths of $G$ cover at least $n-c$ vertices.
\end{conjecture}

 In Section~\ref{nobip}, we prove Conjecture \ref{conj2} for $k=1,t=2$
(Theorem \ref{cexists}) and in particular, $c(C_4,1,2)=1$ (Theorem
\ref{allbutone}). Note that this conjecture is related to Conjecture
\ref{saconj1} by selecting $H$ to be the complete graph of size
$k+1$.

\section{Partitioning into connected
matchings.}\label{matchingresults}

In this section we prove Conjectures \ref{saconj} and \ref{schconj}
in weakened forms, replacing cycles and paths with connected
matchings (Theorems \ref{connmatch}, \ref{degmatch}).
We notice first that the $t=1$ case of Conjecture \ref{saconj} is
due to P\'osa~\cite{posa}.\footnote{See also Exercise 3 on page 63
in \cite{LO}.}

\begin{lemma}\pont\label{posa} The vertex set  of any graph $G$ can be
partitioned into at most $\alpha(G)$ parts, where each part  either
contains  a spanning cycle, or spans an edge or a vertex.
\end{lemma}

For  two colors, we need the following result, which is essentially
equivalent to K\"onig's theorem. It  was discovered in \cite{Gydisz}
and applied in \cite{GYS1}.

\begin{lemma}\pont\label{conncov} Let the edge set of $G$ be colored
with two colors. Then $V(G)$ can be covered with  the vertices of at
most $\alpha(G)$ monochromatic connected subgraphs of $G$.
\end{lemma}

\noindent {\bf Proof.} For a graph $G$ whose edges are colored with
red and blue, let $\rho(G)$ denote the minimum number of
monochromatic components covering the vertex set of $G$. Let
$\alpha^*(G)$ be the maximum number of vertices in $G$ so that no
two of them is covered by a monochromatic component. Suppose that the red edges define connected components $C_1,\dots,C_p$ and the blue edges define connected components $D_1,\dots,D_q$. Define a bipartite multigraph $B$ with vertex classes $C_1,\dots,C_p$ and $D_1,\dots,D_q$. For every vertex $v\in V(G), v\in A_i, v\in B_j$ we define the edge $C_i,D_j$ in $B$.  (In fact, $B$ is the dual of the hypergraph formed by the monochromatic components on $V(G)$.)

Recall that  $\nu(B)$ is the maximum number
of pairwise disjoint edges in $B$ and $\tau(B)$ is the minimum
cover, i.e., the least number of vertices in $B$ that meet all edges of
$B$. From K\"onig's theorem and from easy observations follows that

\begin{equation}\label{connection}
\rho(G)=\tau(B)=\nu(B)=\alpha^*(G)\le\alpha(G)
\end{equation}
finishing the proof. \qed

Observe that (\ref{connection}) gives a stronger form of Lemma \ref{conncov} (equivalent form of K\"onig's theorem).

\begin{proposition}\pont\label{alphastar} For any $2$-edge colored graph $G$, $\rho(G)=\alpha^*(G)$.
\end{proposition}

\begin{theorem}\pont\label{connmatch}
If the edges of a graph $G$ are colored red and blue, then $V(G)$
can be partitioned into at most $2\alpha(G)$ monochromatic parts,
where each part is either an
edge, or a single vertex, or contains a connected matching or a spanning cycle.
\end{theorem}

It is worth noting that Theorem~\ref{connmatch} is best possible,
although it is weaker than Conjecture~\ref{saconj}. Indeed, let $G$
be formed by $k$ vertex disjoint copies of $K_s$, where $s\ge 3$. We color $E(G)$
so that in each $K_s$ the set of blue edges forms a $K_{s-1}$.
Here $\alpha(G)=k$, and we need two parts to cover  each $K_s$, one
in each color.

\noindent {\bf Proof of Theorem~\ref{connmatch}.} Set $V=V(G)$. By
Lemma \ref{conncov}, we can cover $V$ by the vertices of some $p$
red and $q$ blue monochromatic components,
$C_1,\dots,C_p,D_1,\dots,D_q$, where $p+q\le \alpha(G)$. We
partition $V$ into the doubly and singly covered sets. Let
$A_{ij}=C_i\cap D_j$ and $S_i=C_i-\cup_j A_{ij}$, $T_j=D_j-\cup_i
A_{ij}$, where $1\le i \le p, 1\le j \le q$.

Fix  $M_i$, a largest red matching in $C_i$ for every $i$, and then
let $N_j$ be a largest blue matching in $D_j-\cup_i V(M_i)$.
These $p+q\leq \alpha(G)$ monochromatic matchings are connected.
Delete the vertices of  these matchings from $V$ and for convenience
keep the same notation for the truncated sets, so $A_{ij}, S_i, T_j$
denote the sets remaining after all vertices of these matchings are
deleted. Denote the remaining graph by $G_1$, and its vertex set by
$V_1$. Partition $V_1$ into three sets, $A=\cup_{i=1}^p\cup_{j=1}^q
A_{ij}, S=\cup_{i=1}^p S_i, T=\cup_{j=1}^q T_j$. Observe that there are no edges between $S$ and $T$.

Edges of $G_1$ can only be inside  $S$ (colored blue) or
inside  $T$ (colored red).  Applying Lemma \ref{posa} for the
blue and red graphs $G_1[S],G_1[T]$, we can cover $S\cup T$ by
$\alpha(G_1[S])+\alpha(G_1[T])$ parts, where each part contains a
monochromatic spanning cycle or it is  an edge or a vertex. Now $A$ is a
 collection of isolated points in $G_1$; we just cover it with its vertices.  Altogether, we
partitioned $V_1$ into $|A|+\alpha(G_1[S])+\alpha(G_1[T])\le
\alpha(G_1)\le\alpha(G)$ parts and together with the monochromatic
connected matchings $M_i,N_j$, there are at most $2\alpha(G)$ parts
as required. \qed

\begin{theorem}\pont\label{degmatch}
 Let $G=(V,E)$ be an $n$-vertex graph with $\delta(G)\ge 3n/4$, where
$n$ is even. If the edges of $G$ are $2$-colored with red and blue,
then there exist a red connected matching and a vertex-disjoint blue
connected matching, which together form a perfect matching of $G$.
\end{theorem}

\noindent {\bf Proof.}  Let $C_1$ be a largest monochromatic
component, say red. Theorem 1.4 in \cite{GYS2} yields $|C_1|\geq
3n/4$. Let $U=V\setminus V(C_1)$. Any vertex $u$ in $U$ can only
have less than $n/4$ red neighbors. Therefore, the blue degree of
$u$ is at least $n/2$. This implies that  the blue neighborhoods of
any two vertices in $U$ which are not connected with a blue edge
intersect. Therefore, if $U\ne \emptyset$, then $U$ is covered by a
blue component of $G$, say $C_2$. If $U=\emptyset$, then define
$C_2$ as a largest blue component in $G$. Set $p=|V(C_1)\setminus
V(C_2)|, q=|V(C_2)\setminus V(C_1)|$, where $p\ge q$ by the choice
of $C_1$. Let $G_1$ be the graph, which we get from
$G$ by deleting the blue edges inside $C_1\setminus C_2$ and the red edges inside
$C_2\setminus C_1$. Note that in Cases 2 and 3 $C_2\setminus C_1=\emptyset$. We distinguish three cases.

\underline{Case 1:} Suppose $|C_1|<n$. By the maximality of $C_1$
and $C_2$, there are no edges between $C_1\setminus C_2$ and
$C_2\setminus C_1$. Therefore, $q<n/4$ and $p<n/4$. We claim that  $G_1$
satisfies the Dirac-property\footnote{Exercise 21 on page 75 in
\cite{LO}.}, $\delta(G_1)\geq n/2$. Indeed, we deleted at most
$n/4-1$ edges at any vertex, and thus the remaining degree is more
than $n/2$ at each vertex. Therefore, there is a Hamiltonian cycle,
that also contains a perfect matching. This perfect matching
consists of a connected red matching and a connected blue matching
covering $G$.

\underline{Case 2:} Suppose $|C_1|=n$ and $p\le n/2$.
 Now we claim that $G_1$ satisfies the
Chv\'atal-property\footnote{Exercise 21 on page 75 in \cite{LO}.}:
if the degree sequence in $G_1$ is $d_1\leq d_2\leq \ldots \leq
d_n$, then $d_k+d_{n-k}\geq n$ for $k\leq n/2$. Indeed, the degrees
of the $p$ vertices in $C_1\setminus C_2$ are at least $3n/4-p+1$,
where $p\leq n/2$. The rest of the degrees are unchanged being at
least $3n/4$. That yields $3n/4 -p+1+3n/4 =3n/2 -p+1> n$ in the
Chv\'atal-condition. This implies the existence of a Hamiltonian
cycle, which contains a perfect  matching. This perfect matching
contains a connected red matching and a connected blue matching,
which together cover $G$.

\underline{Case 3:} Suppose $|C_1|=n$ and $p>n/2$. That is,
$|C_2|-p<n/2$.  Again, we
claim that there is a perfect matching in $G_1$. Assume to the
contrary that the largest matching is imperfect. By Tutte's theorem,
there exists a set $X$ of vertices in $G_1$ such that the number of
odd components in $G_1\setminus X$ is larger than $|X|$, which  implies that  $|X|<n/2$. Let all the
components (not just the odd ones) be $D_1, D_2, \dots, D_{\ell}$ in
increasing order of their size, $\ell\ge |X|+1$. Note that $\ell\geq 2$
always holds, even for $X=\emptyset$, as $n$ is even.
  Notice,
that any potential edge in $G$ between two components of
$G_1\setminus X$ is a blue edge inside $C_1\setminus C_2$ that was
deleted. Let $H$ be the graph formed by the vertices in $G\setminus
X$, and the blue edges in $C_1\setminus C_2$. Since $|X|<n/2$, we have
$|V(H)|>n/2$.

Suppose first that $|X|=x< n/4$. Let us consider the smallest
component $D_1$ and put $|D_1|=d_1$. We claim that
\begin{equation}\label{n/2}d_1+x\leq n-|D_1\cup X| .\end{equation}
For
$d_1=1$, using $n$ being even, we also get (\ref{n/2}) from
$|D_1\cup X|=1+x\leq n/2$.
When $x=0$ then (\ref{n/2}) is true as $\ell\ge 2$, when $x=1$ then (\ref{n/2}) is true because $n$ is even.
For $x\geq 2$ and  $d_1\geq 2$ we have $|D_1\cup X|=d_1+x\leq d_1 x \leq n-|D_1\cup X|$,
implying (\ref{n/2}).

From (\ref{n/2}), the blue neighborhoods of any two vertices in $D_1$ intersect in
$H$, and $D_1$ is covered by a blue component $C_2'$.
Using $x<n/4$, we get $|C_2'|\geq 3n/4-d_1-x+1+d_1-1 = 3n/4 -x >
n/2$. That is a contradiction since $C_2$ was the largest blue
component and $|C_2|<n/2$.

Now we may assume $n/4\le |X|<n/2$. Since $|X|<n/2$ we have
$V(H)>n/2$. If we prove that $H$ is connected, then we get a
contradiction again, since $C_2$ was the largest blue component, and
$|C_2|<n/2$. Assume to the contrary that we can partition the
vertices of $H$ into $A$ and $B$ with no edges between them. We may
assume $|A|\ge |B|$, and therefore $|A|> n/4$. We have two subcases.

\underline{Case 3.a:} Suppose $A\cap D_i\neq \emptyset$ for $1\le
i\le \ell$. Let $v$ be a vertex in $B$ and assume $v\in D_j$. There
is no edge of $G$  from $v$ to $A\cap D_i$, for each $i\neq j$, $1\le i\le
\ell$: An edge from $G_1$ is impossible, because $i\neq j$; a blue
edge from $C_1\setminus C_2$ is impossible, because $(A,B)$ is a
cut in $H$. Therefore, the degree of $v$ in $G$ is at most $n-1-{\ell}+1\le
n-(|X|+1)\le n-1-n/4<3n/4$, a contradiction.

\underline{Case 3.b:} Suppose $A\cap D_j=\emptyset$ for a fixed $j$,
$1\le j\le \ell$. Let $v$ be a vertex in $D_j$. There is no edge
from $v$ to any vertex $u$ of $A$: An edge from $G_1$ is impossible,
because $u\in D_i$, where $i\neq j$. A blue edge from $C_1\setminus
C_2$ is impossible, because $(A,B)$ is a cut. Therefore, the degree
of $v$ in $G$ is at most $n-1-|A|\le n-1-n/4<3n/4$, a contradiction.
\qed

\section{Applying the Regularity lemma.}\label{reg}

As in many applications of the Regularity Lemma, one has to handle
irregular pairs, that translates to exceptional edges in the reduced
graph. To prove  such a variant of Theorem \ref{connmatch}, first
Lemma \ref{conncov} is tuned up. A graph $G$ on $n$ vertices is
$\varepsilon$-{\it perturbed} if at most $\varepsilon {n\choose 2}$
of its edges are marked as exceptional (or perturbed). For a perturbed graph $G$,
let $G^-$ denote the graph obtained by removing all perturbed edges.

\begin{lemma}\pont\label{conncovexp}
Suppose that $G$ is a $2$-edge-colored $\varepsilon$-perturbed graph
on $n$ vertices, $n\ge \varepsilon^{-1/2}$. Then all but at most
$f(\alpha(G))\sqrt{\varepsilon}n$ vertices of $G$  can be covered by
the vertices of  $\alpha(G)$ monochromatic connected subgraphs of
$G^-$, where $f$ is a suitable function.
\end{lemma}

\noindent {\bf Proof.} Set $\alpha=\alpha(G)$ and remove from $V(G)$
a set $X$ of at most $\sqrt{\varepsilon}n$ vertices so that in the
remaining graph $H$ each vertex is incident to at most
$\sqrt{\varepsilon}n$ perturbed edges.

Let $\cal{T}$ denote the (possibly edgeless) hypergraph whose edges
are those sets  $T\subset V(H)$ for which $|T|=\alpha+1$ and no
monochromatic component of $H^-$ covers more than one vertex of $T$. (Each $T\in \cal{T}$ is a witness showing $\alpha^*(H^-)\ge \alpha+1$.)
We call pairwise disjoint hyperedges $T_1,T_2,\dots,T_k$ in
$\cal{T}$ {\it independent}, if there are no perturbed edges in the $k$-partite graph defined by the $T_i$-s. Set
$c=3^{\alpha^2}$ and let $R=R(3,3,\dots,3,\alpha+1)$ be the $c$-color Ramsey number, the smallest $m$ such that in every $c$-coloring of the edges of $K_m$ either there is a triangle in one of the first $c-1$ colors or a $K_{\alpha+1}$ in color $c$.

\begin{claim}\pont\label{ramsey} Select in $\cal{T}$
as many pairwise independent hyperedges as possible, say $T_1,T_2,
\dots,T_k$. Then $k < R$.
\end{claim}

{\bf Proof.} Fix an ordering within each of the sets $T_i$; if $x\in
T_i$ is the $j$-th element in this order in $T_i$, we write
$ind(x)=j$. Suppose for contradiction that $k\ge R$ and consider a
coloring of the pairs among $T_1,T_2,\dots,T_k$ defined as follows.
Color a pair $T_i,T_j$ ($1\le i<j \le k$) by their ``color
pattern'' on the pairs $x\in T_i,y\in T_j$ with $ind(x)\ne ind(y)$. There are $\alpha^2$ such pairs (none of them is a perturbed edge) thus $x,y$ is a red edge, a blue edge or not an edge in $H$. So we have a $c$-coloring on the pairs $T_i,T_j$, the color when all the $\alpha^2$ pairs are not edges of $H$ is called special. By the
assumption $k\ge R$, we have either $\alpha+1$ $T_i$-s
with any pair of them colored with the special color or three $T_i$-s with all three pairs colored with the same non-special color. We show that both cases lead to contradiction.

In the latter case we have a triple, say $T_1,T_2,T_3$ and different indices $i,j$, such that $p\in T_1,q,r\in T_2,s\in T_3$,  $ind(p)=ind(q)=i,ind(r)=ind(s)=j$ and $pr,ps,qs$ are all edges of $H$ colored with the same color. Thus $r,p,s,q$  is a monochromatic path of $H^-$, intersecting $T_2$ in two vertices, contradicting to the definition of $T_2$.

In the former case we have say $T_1,T_2,\dots,T_{\alpha+1}$ pairwise colored with the special color. For
$i=1,2,\dots,\alpha+1$, select $v_i\in T_i$ such that $ind(v_i)=i$. Observe that $\{v_1,\dots,v_{\alpha+1}\}$ spans an independent set
in $G$, contradicting the assumption that $\alpha(G)=\alpha$. $\qed$\\

Let $Y$ denote the set of vertices in $H$ sending at least one perturbed edge to
$\cup_{i=1}^k T_i$. Observe that $|Y|\le (\alpha+1)R\sqrt{\varepsilon} n$ and by the maximality of  $k$,
$Z=\cup_{i=1}^k T_i \cup Y$ meets all edges of $\cal{T}$, thus
removing $X\cup Z$ from V(G) leaves a subgraph $F\subset G$ with
$\alpha^*(F^-)\le \alpha$. Therefore, applying Proposition \ref{alphastar} to $F^-$, $\rho(F^-)\le \alpha$. The theorem
follows, since (using the assumption $1\le \sqrt{\varepsilon}n$)
$$|X\cup Z|\le \sqrt{\varepsilon}n+R(\alpha+1)+(\alpha+1)R\sqrt{\varepsilon}n\le (1+2R(\alpha+1)\sqrt{\varepsilon}n,$$
i.e.
$f(\alpha)=(1+2R(\alpha+1))$ is a suitable function. \qed\\



Now we are ready to prove a {\em perturbed} version of
Theorem~\ref{connmatch}.

\begin{theorem}\pont\label{connmatchpert} Let  $G$ be
an $\varepsilon$-perturbed $2$-edge-colored graph on $n$ vertices, $n\ge \varepsilon^{-1/2}$.
Then there exists a $Z\subset V(G)$ such that $|Z|\le
(f(\alpha(G))+\alpha(G))\sqrt{\varepsilon}n$ and $V(G)\setminus Z$
can be partitioned into at most $2\alpha(G)$ classes, where each
part in $G^-$ either contains  a connected monochromatic spanning  matching or
a monochromatic spanning cycle or it is  an edge or a single vertex.
\end{theorem}

\noindent {\bf Proof.} Using Lemma \ref{conncovexp}, we can remove
from $V(G)$ a set of at most $f(\alpha)\sqrt{\varepsilon}n$ vertices
such that for the remaining graph $H$, the following holds. The
vertices $V(H)$ can be covered by the vertices of at most
$\alpha(G)$ monochromatic components of $H^-$, say  with  $p$ red and
$q$ blue monochromatic components, $C_1,\dots,C_p,D_1,\dots,D_q$,
where $p+q\le \alpha(G)$. We may suppose that each vertex of $H$ is
incident to at most $\sqrt{\varepsilon}n$ perturbed edges, as this
is automatic from the proof of Lemma \ref{conncovexp}. The $p+q$
components yield a partition of $V(H)$ into doubly and singly covered sets. Let
$A_{ij}=C_i\cap D_j$ and $S_i=C_i-\cup_j A_{ij}$, $T_j=D_j-\cup_i
A_{ij}$, where $1\le i \le p, 1\le j \le q$. First let $M_i$ be a
largest red matching induced by $H^-$ in $C_i$ for every $1\le i \le
p,$ and then $N_j$ be a largest blue matching induced by $H^-$ in
$D_j-\cup_i V(M_i)$, for every $ 1\le j \le q$.
 Observe that these matchings are connected in
$H^-$. Delete all vertices of these matchings from $V(H)$ and for
convenience keep the same notation for the truncated sets (so
$A_{ij}, S_i, T_j$ denotes the sets remaining after all vertices of
these matchings are deleted). The remaining graph is denoted by $F$.
Partition $V(F)$ into three sets, $A=\cup_{i=1}^p\cup_{j=1}^q
A_{ij}, S=\cup_{i=1}^p S_i, T=\cup_{j=1}^q T_j$. Observe that edges
of $F^-$ can be only inside $S$ (colored blue) or inside $T$
(colored red). Now we follow the proof method of Lemma \ref{posa}
(see Exercise 3 on page 63 in \cite{LO}) to partition most of the
vertices in $V(F)$ into at most $\alpha(G)$ monochromatic cycles.

We apply the following procedure to subsets $U$ of one of the sets
$A,S,T$. Observe that $F^-[U]$ is an independent set if $U\subset
A$, edges of $F^-[U]$ are all blue if $U\subset S$, edges of
$F^-[U]$ are all red if $U\subset T$.

In any step of the procedure, consider a maximal path $P$ of
$F^-[U]$ and let $x$ be one of its endpoints. If $x$ is an isolated
vertex in $F^-[U]$, define $C^*=\{x\}$. If $x$ has degree one in
$F^-$, let $y$ be its neighbor on $P$ and define $C^*=\{x,y\}$. If
$x$ has degree at least two in $F^-$, let $z$ be the neighbor of $x$
on $P$ (in $F^-$), which is the furthest from $x$. Now $C^*$ is
defined as the cycle obtained by connecting the endpoints of the
edge $xz$ on the path $P$. Let $Y$ be the set of perturbed
neighbors of $x$ in $F^-$. That is, the set of vertices in $V(F)$,
which are adjacent to $x$ by exceptional edges. The step ends with
removing $C^*\cup Y$ from $V(F)$ and defining the new $F,A,S,T$ as
the truncated sets.

This procedure decreases $\alpha(F)$ at each step, because any
independent set of the truncated set can be extended by $x$ to an
independent set of $F$. Therefore, at most $\alpha(G)$ steps can be
executed. Now apart from the union of the sets $Y$s, at most
$\alpha(G)$ monochromatic $C^*$-s partition $V(F)$. Together with
the $p+q\le \alpha$ monochromatic connected matchings $N_i,M_j$ we
have the required covering. The number of uncovered vertices are at
most $f(\alpha)\sqrt{\varepsilon}n$ (lost when the matchings were
defined) plus $\alpha \sqrt{\varepsilon}n$ (when the cycles are
defined). \qed

\subsection{Building cycles from connected matchings.}\label{reglemma}

Next we show how to prove Theorem~\ref{saass} from Theorem
\ref{connmatchpert} and the Szemer\'edi Regularity Lemma~\cite{Sz}. The material
of this section is fairly standard by now (see \cite{rlogr,GYLSS,GRSSz,GYRSSZ,GRSSZ}
so we omit some of the
details. We need a $2$-edge-colored version of the Szemer\'edi Regularity
Lemma.\footnote{For background, this variant and other variants of
the Regularity Lemma see \cite{KS}.}

\begin{lemma}\label{2-reg}\pont
For every integer $m_0$ and positive $\varepsilon$, there is an
$M_0=M_0(\eps, m_0)$ such that for $n\geq M_0$ the following holds.
For any $n$-vertex  graph $G$, where $G=G_1\cup G_2$ with
$V(G_1)=V(G_2)=V$, there is a partition of $V$ into $\ell+1$
clusters $V_0,V_1,\dots,V_\ell$ such that
\begin{itemize}
\item $m_0\leq \ell\leq M_0$, $|V_1|=|V_2|=\dots=|V_\ell|$,
$|V_0|<\eps n$,
 \item apart from at most $\eps {\ell
\choose 2}$ exceptional pairs, all pairs $G_s|_{V_i\times V_j}$ are
$\eps$-regular, where $1\leq i<j\leq \ell$ and $1\leq s\leq 2$.
\end{itemize}
\end{lemma}

{\bf Proof of Theorem \ref{saass}.} Given $\eta$ and $\alpha$, first
we fix a  positive $\varepsilon$ sufficiently small so that the
claimed bound $(f(\alpha)+\alpha)\sqrt\varepsilon$ in
Theorem~\ref{connmatchpert} is much smaller than $\ \eta$. Then we
choose $m_0$ sufficiently large compared to $1/\sqrt\varepsilon$ (so
Theorem~\ref{connmatchpert} can be applied). Let $G$ be a graph on
$n$ vertices with $\alpha(G)=\alpha$, where $n\geq M_0$ with $M_0$
coming from Lemma~\ref{2-reg}.
 Consider a
$2$-edge-coloring of $G$, that is $G=G_1\cup G_2$. We apply
Lemma~\ref{2-reg} to $G$ in order to obtain a partition of
$V$, that is $V=\cup_{0\leq i\leq \ell}V_i$. Define the following
{\em reduced graph} $G^R$: The vertices of $G^R$ are
$p_1,\ldots,p_\ell$, and there is an edge between vertices $p_i$ and
$p_j$ if the pair $(V_i, V_j)$ is either exceptional\footnote{That
is, $\eps$-irregular in $G_1$ or in $G_2$. Also,   these edges are marked
exceptional in $G^R$.}, or if it is $\eps$-regular in both $G_1$ and $G_2$
with density in $G$ exceeding $1/2$. The edge $p_ip_j$ is colored
with the color, which is used on the most edges from $G[V_i, V_j]$ (the
bipartite subgraph of $G$ with edges between $V_i$ and $V_j$). The
density of this majority color is still at least $1/4$ in $G[V_i,
V_j]$. This defines a $2$-edge-coloring $G^R=G_1^R\cup G_2^R$.

We claim that $\alpha(G^R)\leq \alpha(G)=\alpha$.
Indeed, we apply the standard Key Lemma\footnote{Theorem~2.1 in
\cite{KS}.} in the complement of $G^R$ and $G$. Note that a
non-exceptional pair is $2\eps$-regular in $\overline{G}$ as well. If
we had an independent set of size $\alpha+1$ in $G^R$, then we would
have an independent set of size $\alpha+1$ in $G$, a contradiction.

We now apply Theorem~\ref{connmatchpert} to the $\eps$-perturbed
$2$-edge-colored $G^R$ (note that the condition in
Theorem~\ref{connmatchpert} is satisfied since $\ell \gg
1/\sqrt\varepsilon$). We cover most of $G^R$ by at most
$2\alpha(G^R)\leq 2\alpha(G)=2\alpha$ subgraphs of $(G{^R})^-$,
where each subgraph in $(G{^R})^-$ is either a connected
monochromatic matching or a monochromatic cycle or an edge or a
single vertex. Finally, we lift the connected matchings back to
cycles in the original graph using the following\footnote{As in
\cite{GYLSS,GRSSz,GYRSSZ,GRSSZ}.} lemma in our context, completing
the proof. Indeed, the number of vertices left uncovered in $G$ is
at most $$(f(\alpha)+\alpha)\sqrt\varepsilon n + 3 \epsilon n +
\epsilon n = (f(\alpha)+\alpha)\sqrt\varepsilon n + 4 \epsilon n
\leq \eta n,$$ using our choice of $\epsilon$. Here the uncovered
parts come from Theorem~\ref{connmatchpert}, from Lemma \ref{lift}
and $V_0$. \qed

\begin{lemma}\pont\label{lift}
Assume that there is  a monochromatic connected matching $M$ (say in
$(G_1{^R})^-$) saturating at least $c|V(G{^R})|$ vertices of
$G{^R}$, for some positive constant $c$. Then in the original $G$
there is  a monochromatic cycle in $G_1$ covering at least
$c(1-3\eps)n$ vertices.
\end{lemma}

{\bf Proof of Theorem \ref{mindegass}.} We combine the degree form
and the $2$-edge-colored version of the Regularity Lemma.

\begin{lemma}\label{dreg}\pont
For every positive $\eps$ and integer $m_0$, there is an
$M_0=M_0(\eps, m_0)$ such that for $n\geq M_0$ the following holds.
For any $n$-vertex graph $G$, where $G=G_1\cup G_2$ with
$V(G_1)=V(G_2)=V$, and real number $\rho \in[0,1]$, there is a
partition of $V$ into $\ell+1$ clusters $V_0,V_1\dots,V_\ell$, and
there are subgraphs $G'=G_1'\cup G_2'$, $G_1'\subset G_1$,
$G_2'\subset G_2$
with the following properties:
\begin{itemize}
\item $m_0\leq \ell\leq M_0$, $|V_0|\leq\eps|V|$,
$|V_1|=\dots=|V_\ell|=L$, \item $deg_{G'}(v)>deg_G(v)-(\rho
+\eps)|V|\text{for all}v\in V$,
\item the vertex sets $V_i$ are independent in $G'$,
\item
each pair $G'|_{V_i\times V_j}$ is $\eps$-regular, $1\leq i<j\leq
\ell$, with  density $0$ or exceeding $\rho$,
\item
each pair $G_s'|_{V_i\times V_j}$ is $\eps$-regular, $1\leq i<j\leq
\ell, 1\leq s\leq 2$.
\end{itemize}
\end{lemma}

Let  $\eps\ll \rho \ll \eta \ll 1$, $m_0$ sufficiently large
compared to $1/\varepsilon$ and $M_0$ obtained from
Lemma~\ref{dreg}. Let $G$ be a graph on $n>M_0$ vertices with
$\delta(G)>(\frac{3}{4}+\eta) n$. Consider a $2$-edge-coloring of
$G$, that is $G=G_1\cup G_2$. We apply Lemma~\ref{dreg} to $G$. We
obtain a partition of $V$, that is $V=\cup_{0\leq i\leq \ell}V_i$.
We define the following {\em reduced graph} $G^R$: The vertices of
$G^R$ are $p_1, \ldots , p_\ell$, and there is an edge between
vertices $p_i$ and $p_j$ if the pair $(V_i, V_j)$ is $\eps$-regular
in $G'$ with density exceeding $\rho$. Since $\delta(G')
> (\frac{3}{4} +\eta  - (\rho + \eps)) |V|$,
calculation\footnote{See a similar computation in \cite{SG}.} shows
that $\delta(G^R) \geq \left(\frac{3}{4} +\eta - 2\rho \right)\ell>
\frac{3}{4}\ell$. The edge $p_ip_j$ is colored again with the
majority color, and the density of this color is still at least
$\rho/2$ in $K(V_i, V_j)$.

Applying Theorem \ref{degmatch} to $G^R$, we get a red connected
matching and a vertex-disjoint blue connected matching, which
together form a perfect matching of $G^R$. Finally we lift the
connected matchings back to cycles in the original graph using
Lemma~\ref{lift}. The number of vertices left uncovered in $G$ is at
most $\sqrt\varepsilon n \leq \eta n$. \qed

\section{Excluding bipartite graphs from the
complement.}\label{nobip}

In what follows, we prove the $t=2,k=1$ case of
Conjecture~\ref{conj2}. As every bipartite graph is a subgraph of a
complete bipartite graph, we may assume that the graph $H$ forbidden
in the complement of $G$ is $K_{p,p}$. Note that the constant $c$ we
get could be greatly improved even using the same arguments with
more involved calculations, however, it would be still far from
being optimal. We use the following well-known theorems.

\begin{theorem}[Erd\H os-Gallai \cite{EG}]\pont\label{erdga}\footnote{See also Exercise 28 on page 76 in
\cite{LO}.} If $G$ is a graph on $n$ vertices with $|E(G)|>
\ell (n-1)/2$, then $G$ contains a cycle of length at least $\ell+1$.
\end{theorem}

\begin{theorem}[K\H ov\'ari-T. S\'os-Tur\'an \cite{KOV}]\pont\label{kotstu}\footnote{See also Exercise 37 on page 77 in
\cite{LO}.} If $G$ is a graph on $n$ vertices such that $K_{p,p}$ is
not a subgraph of $G$, then $|E(G)| \le
(p-1)^{1/p}n^{2-1/p}+(p-1)n\le 2p n^{2-1/p}$.
\end{theorem}

\begin{lemma}\pont\label{twopossib}
Let $p$ and $n$ be positive integers such that $n\ge (10p)^p$. Let
$G$ be an $n$-vertex graph such that $K_{p,p} \not \subset
\overline{G}$. Then any $2$-edge-coloring of $G$ contains a
monochromatic cycle of length at least $n/4$.
\end{lemma}

\noindent {\bf Proof.} By Theorem~\ref{kotstu} and by the lower
bound on $n$,  $$e(G)\ge {n\choose 2}-2pn^{2-1/p}=
n^2/2-n/2-2pn^{2-1/p}\ge n^2/2-n/2-n^2/5 \ge n^2/4,$$ so one of the
colors, say red, is used at least $n^2/8$ times.
Then using Theorem \ref{erdga} in the red subgraph we get a red cycle of length at least $n/4.$\qed\\

For a bipartite graph $G$ with classes $A,B$, the {\it bipartite
complement} $\overline G[A,B]$ of $G$ is obtained via complementing
the edges between $A$ and $B$, and keeping $A$ and $B$ independent
sets.

\begin{lemma}\pont\label{bipart} Let $0<\epsilon<1$ and $n \ge (50p)^p/\epsilon$.
Let $G$ be a bipartite graph with classes $A$ and $B$, $|A|=|B|=n$
such that $ K_{p,p} \not\subset \overline G[A,B]$.  Then there is a
path of length at least $(2-\epsilon)n$ in $G$.
\end{lemma}

\noindent {\bf Proof.} First we prove a weaker statement.

\begin{claim}\pont\label{pathh} Let $G'$ be a bipartite graph with classes
$A'$ and $B'$ with $|A'|=|B'|=m\ge (20p)^p$ such that $ K_{p,p}
\not\subset \overline G'[A',B']$.  Then there is a path of length at
least $m/2$ in $G'$.
\end{claim}
{\bf Proof.} By Theorem~\ref{kotstu}, $e(G')\ge
m^2-8pm^{2-1/p}>m^2/2=(2m)^2/8$, so by Theorem \ref{erdga} $G'$
contains a path of length at least $m/2$.\qed\\

Let $P$ be a longest path in $G$. Using  Claim~\ref{pathh} with
$G=G'$, we have that $|P|\ge n/2$. Assume for a contradiction that
$P$ is shorter than $(2-\epsilon)n$. Because $G$ is bipartite, we
can choose $A'\subset (G-P)\cap A$ and $B'\subset (G-P)\cap B$ with
$|A'|=|B'|> \epsilon n/3$. By Claim~\ref{pathh},
 $G[A',B']$ contains a path $P'$ with at least $\epsilon n/6$ vertices.

Consider the last $2p$ vertices of $P$ and the last $2p$ vertices of
$P'$. There is an edge $e$ between these set of vertices by the
assumption. Adding $e$ to $P\cup P'$, there is a path, which
contains all but $2p$ vertices of $P$, and all but $2p$ vertices of
$P'$,  hence it is longer than $P$, a contradiction. Here we used
that $\epsilon n/6> 4p$. \qed

\begin{theorem}\pont\label{cexists}
Let $G$ be an $n$-vertex graph such that $K_{p,p}\not \subseteq
\overline{G}$. Then any $2$-edge-coloring of $G$ contains two vertex
disjoint monochromatic paths of distinct colors covering at least
$n- 1000(50p)^p$ vertices.
\end{theorem}

\noindent {\bf Proof.} Consider the vertex disjoint blue path, red
path pair $(P_1,P_2)$, which cover the most vertices, and let
$G'=G\setminus \{P_1\cup P_2\}$. Suppose there are $n_1$ vertices in
$G'$, where $n_1>1000(50p)^p$. As $n>n_1>1000(50p)^p$,
 by  Lemma~\ref{twopossib} at least $n/4$ vertices are covered by $P_1\cup P_2$.
Let $t=10(50p)^p<n_1/100.$ We split the proof into two cases.

\underline{Case 1:} One of the paths, $P_2$ say, is shorter than
$t$. Using that $3t< n/4$ we have that the length of $P_1$ is at
least $2t$ in this case. Now $G'$ does not contain a red path of
length $t$, but by Lemma~\ref{twopossib} it contains a monochromatic
cycle of length at least $n_1/4>4t$, which must be blue. Hence, $G'$
contains a blue path, say  $P_3$, of length at least $4t$.

Denote $L_1$, the set of last $2t$ vertices of $P_1$ and $L_3$, the
set of  last $2t$ vertices of $P_3$. There is an edge $e$ between
$L_1$ and $L_3$ as $2t >p$ and $K_{p,p}\not \subseteq \overline{G}$.
If $e$ was blue then we use $e$ to connect the paths $P_1,P_3$, and we find
a blue path longer than $P_1$ vertex disjoint from $P_2$, a
contradiction.

Hence all edges between $L_1$ and $L_3$ are red, and we can apply
Lemma \ref{bipart} for the red bipartite graph between $L_1$ and
$L_3$ with $\epsilon=1/8$. (Note that $2t\geq 8 (50p)^p$, so indeed
the lemma is applicable.) It yields a red path $P_4$ of length
$(2-1/8)2t$ in $L_1\cup L_3$. Let $P_1'$ be $P_1$ without the last
$2t$ vertices. Now $P_1'$ and $P_4$ are disjoint and cover
 more vertices than $P_1$ and $P_2$, which is a contradiction.

\underline{Case 2:} Both $P_1$ and $P_2$ have length at least $t$.
Without loss of generality, in $G'$ Lemma~\ref{twopossib} implies
the existence of a blue cycle $C$ of length at least $n_1/4\ge 4t$.
Denote $R_1$ the set of the last $t$ vertices of $P_1$, $R_2$ the
set of the last $t/2$ vertices of $P_2$, and $C_1$ any set of
consecutive $t$ vertices of $C$. There are no blue edges between
$R_1$ and $C_1$, otherwise $P_1$ could be replaced with a longer
blue path. Now by Lemma \ref{bipart}, with $\epsilon=1/8$, there is
a red path $P_3$ in $G(R_1,C)$ of length $15t/8$. Let $B$ be the set
of the first and last $t/4$ vertices of $P_3$. For each vertex $v$
in $B$, there is a red path $P_v$ of length $13t/8$ starting at $v$,
which is a subpath of $P_3$. If there is a red edge $e=(u,v)$
between $R_2$ and $B$, then $P_2\cup e\cup P_v$ contains a red path
with at least $|P_2|+13t/8-t/2$ vertices which together with the
disjoint $P_1-R_1$ cover more vertices than the pair $(P_1,P_2)$, a
contradiction.

Therefore, there are only blue edges between $B$ and $R_2$. Since
$|B\cap P_1|\ge p$, there are at least $t/2-p+1$ vertices of $R_2$
having neighbors in $B \cap P_1$. Let $R_2'$ be the set of those
vertices. If there is a blue edge $f$ between $R_2'$ and $C$, then
$P_1\cup f\cup C$ contains a blue path which together with the
disjoint $P_2-R_2$ cover more vertices than the pair $(P_1,P_2)$, a
contradiction.

Therefore, all the edges between $R_2'$ and $C$ are red. We already
know that there are no red edges from $R_2'$ and $B\cap C$. But we
have that $|R_2'|\geq p$ and $|B\cap C|\geq p$, which is a
contradiction.
\qed\\

The following proposition, which is a $1$-colored version of one of our main results, Theorem~\ref{allbutone}, is also a special case of $R(P_m,C_n)$,
determined in \cite{FLPS}.

\begin{proposition}\pont\label{pathc4ramsey}
If $G$ is a graph on $n$ vertices and $C_4\not\subseteq
\overline{G}$, then $G$ contains a path, which covers $n-1$
vertices.
\end{proposition}

\noindent {\bf Proof.}  Denote by $P$ a  longest path  of $G$. Let
$a$ and $b$ be the first and last vertex of $P$. If $P$ contains
less than $n-1$ vertices, then there are two vertices $x$ and $y$
not in $P$. Let us consider the pairs $ax$, $xb$, $by$, $ya$. If
none of them spans an edge  in $G$, then they span a $C_4$ in
$\overline{G}$, which is a contradiction. If any of them spans an
edge in $G$, then it extends $P$,
which is again a contradiction. \qed\\

The following result, the two-color version of Proposition
\ref{pathc4ramsey}, shows that Conjecture \ref{conj2} is true for
$H=C_4$ with $c(C_4)$=1.

\begin{theorem}\pont\label{allbutone}
Let $G$ be a graph such that $|V(G)|\ge 7$ and $C_4\not\subseteq
\overline{G}$. If the edges of $G$ are colored red and blue, then
there exist two vertex-disjoint monochromatic paths of different
colors covering $n-1$ vertices.
\end{theorem}

For simplicity, we refer to edges of $\overline{G}$ as black edges,
and think of $G$ as $K_n$ with a $3$-edge-coloring, but
monochromatic paths should be blue or red, and sometimes when we
write "edge of G" we mean "red or blue edge of $G$". We trust that
this will not confuse the reader.

{\bf Remark 2.} The value $n-1$ in Theorem \ref{allbutone} is best
possible, as shown by the following example. Let $v_1$ and $v_2$ be
two different vertices in $K_n$. If $v_1x$ is black for all $x$, and
$v_2y$ is red for all $y$, $y\in V(K_n)\setminus v_1$, and all other
edges are blue, then any two monochromatic paths can only cover at
most $n-1$ vertices.

The condition $|V(G)|\ge 7$ is somewhat unexpected, since the
statement is true if $|V(G)|\le 4$. On five vertices, let
$G_5=K_1\cup C_4$ and color the edges of $C_4$ alternately red and
blue. On six vertices, let $G_6$ be the complement of $C_6$ and
color the long diagonals red and the short diagonals blue. One can
easily check that pairs of vertex disjoint red and blue paths must
leave two vertices uncovered in these
graphs. \\

\noindent {\bf Proof Theorem \ref{allbutone}.} Fix a blue path
$P_1=a_1 \dots a_i$ and a red path $P_2=b_1 \dots b_j$ such that
$i+j$ is as large as possible, and under this condition  $|i-j|$ is
as small as possible. Let $G'$ be $G \setminus (P_1 \cup P_2)$. If
$G'$ contains only one vertex, then we are done. Therefore, we may
choose a  $U\subseteq V(G')$ such that  $U=\{x,y\}$ for some $x \neq
y$. Since $i+j$ is maximal, there are no blue edges between
$\{a_1,a_i\}$ and $G'$ and there are no red edges between
$\{b_1,b_j\}$ and $G'$. We consider two cases, according whether
$\min{i,j}=1$ (say then $i=1$).

\underline{Case 1:} $i=1$. If there is a blue edge between $b_1$ and
$G'$, then that one edge and $b_2 \dots b_j$ would be a better pair
of paths (with smaller difference of the sizes), which is a
contradiction, unless $j=2$. In this case, $X=V(G)\setminus
\{b_1,b_2\}$ has at least five vertices and (using that no $C_4$ in
$\overline G$) one can easily see that $X$ has either a blue edge or
a red $P_3$ and both contradicts the choice of $P_1,P_2$.

\underline{Case 2:} $i,j \ge 2$. Since there is no black $C_4$,
there is an non-black edge of $G$ between some of the endpoints of
$P_1$ and some of the  endpoints of $P_2$. We call such an edge  a
{\it cross-edge}.

\begin{claim}\pont\label{ends}
If both endpoints of a cross-edge are connected to $G'$ by a
non-black edges of $G$, then we can increase the number of vertices
covered by the two monochromatic paths.
\end{claim}

We may assume that   $a_1b_1$ is a cross-edge and it is blue. There
is a blue edge between $b_1$ and $G'$, say $b_1z$. Now $zb_1a_1\dots
a_i$ and $b_2\dots b_j$ are two monochromatic paths,
which cover more vertices than $P_1$ and $P_2$. $\Box$\\


In what follows, we may assume that $a_1b_1$ is a blue cross-edge,
and $b_1z$ is black for any vertex $z$ of $G'$. Let $v\in V(P_1)
\cup V(P_2)\setminus b_1$. If $vz_1$ and $vz_2$ were two black edges
for some  $z_1,z_2\in G'$, then $vz_1b_1z_2$ would be a black
$4$-cycle, a contradiction. Therefore, $v$ is adjacent to all but
one vertex in $G'$. In particular, there are red edge from both
$a_1$ and $a_i$ to $G'$ and a blue edge from $b_j$ to $G'$.
Therefore, the edges $a_1b_j$ and $a_ib_j$ are both black by
Claim~\ref{ends}.

\underline{Case 2.1:} $j=2$. If there were two (red) edges between
$a_i$ and $G'$, say $a_iz_1$ and $a_iz_2$, then $b_1a_1\dots
a_{i-1}$ and $z_1a_iz_2$ would cover more vertices than $P_1\cup
P_2$, a contradiction. Therefore, $|V(G')|=2$, that is $U=G'$. We
may assume $a_ix$ is red and $a_iy$ is black. It follows that $a_1y$
is red and $a_1x$ is black, otherwise
$a_1ya_ib_j$ would be a black $C_4$. 
contradiction. Since $|V(G)|\ge 7$, we now get $i>2$. Therefore,
$a_{i-1}\ne a_1$.

\underline{Case 2.1.1:} $a_1x$ is black. Consider the edges
$a_{i-1}x$ and $a_{i-1}b_2$. If both of them were black, then
$a_1xa_{i-1}b_2$ would be a black
$C_4$. 
If both of them were red, then $b_1b_2a_{i-1}xa_i$ and $a_1\dots
a_{i-2}$ would cover more vertices than $P_1\cup P_2$. If
$b_2a_{i-1}$ is blue, then $b_1a_1\dots a_{i-1}b_2$ and $a_ix$ cover
more vertices than $P_1\cup P_2$.

If $a_{i-1}x$ is blue, then consider the existing blue edge between
$b_2$ and $U$. If $b_2x$ were blue, then $b_1a_1\dots a_{i-1}xb_2$
and $a_i$ would cover more vertices than $P_1\cup P_2$. Therefore,
$b_2y$ is a blue edge. Consider now the edge $b_1a_i$. If $b_1a_i$
were red, then $b_2b_1a_ix$ and $a_1\dots a_{i-1}$ would cover more
vertices than $P_1\cup P_2$. If $b_1a_i$ were blue, then
$xa_{i-1}a_ib_1a_1\dots a_{i-2}$ and $b_2$ would cover more vertices
than $P_1\cup P_2$. Therefore, $b_1a_i\in \overline G$. Now we
consider the edge $xy$. If $xy$ is blue, then $a_1\dots a_{i-1}xy$
and $b_1b_2$ cover more vertices than $P_1\cup P_2$. If $xy$ is red,
then $b_1a_1\dots a_{i-1}$ and $a_ixy$ cover more vertices than
$P_1\cup P_2$. Finally, if $xy\in \overline{G}$, then $xya_ib_1$ is
a black 4-cycle.
This shows that $a_{i-1}x$ is not blue.

Now one of $a_{i-1}x$ and $a_{i-1}b_2$ is red and the other one is
black. If $a_{i-1}x$ is red, then consider  $a_{i-1}y$. If
$a_{i-1}y$ is red, then $a_ixa_{i-1}y$ and $b_1a_1\dots a_{i-2}$
cover more vertices than $P_1\cup P_2$. If $a_{i-1}y$ is blue, then
$b_1a_1\dots a_{i-1}y$ and $a_ix$ cover more vertices than $P_1\cup
P_2$.
If $a_{i-1}y\in \overline{G}$, then $b_2a_{i-1}ya_i$ is a black $4$-cycle.\\
If $a_{i-1}b_2$ is red and $a_{i-1}x$ is black, then look at
$a_{i-1}y$. If $a_{i-1}y$ is black, then $xa_{i-1}yb_1$ is a black
$C_4$. If $a_{i-1}y$ is blue, then $b_1a_1\dots a_{i-1}y$ and $a_ix$
cover more vertices than $P_1\cup P_2$. If $a_{i-1}y$ is red, then
$b_1a_1\dots a_{i-2}$ and $b_2a_{i-1}y$ cover the same number of
vertices as $P_1\cup P_2$. At the same time, if $i\ge 4$, $|i-j|$ is
smaller, giving a contradiction. On the other hand, if $i=3$, then
$a_i$ and $b_1b_2a_{i-1}ya_1$ and $a_i$ cover more vertices than
$P_1\cup P_2$.

\underline{Case 2.1.2:} $a_1x$ is red. If $i\ge 4$, then $a_2\dots
a_{i}$ and $xa_1y$ cover the same number of vertices as $P_1\cup
P_2$ with a smaller $|i-j|$, a contradiction. Therefore, $i=3$ that
is $|V(G)|=7$. If $b_2y$ is blue, then look at $a_2y$. If $a_2y$ is
blue, then $b_1a_1a_2yb_2$ and $a_3x$ cover more vertices than
$P_1\cup P_2$. If $a_2y$ is red, then $b_1$ and $a_3xa_1ya_2$ cover
more vertices than $P_1\cup P_2$. Therefore, $a_2y$ is black. Now if
$a_2b_2$ is black, then $b_2a_3ya_2$ is a black $C_4$. If $a_2b_2$
is blue, then $b_1a_1a_2b_2y$ and $a_3x$ cover more vertices than
$P_1\cup P_2$. Therefore, $a_2b_2$ is red. Now $a_2x$ must be blue
and $b_2x$ black. Consider now $b_1a_3$. If $b_1a_3$ is blue, then
$a_3b_1a_1a_2x$ and $b_2$ cover more vertices than $P_1\cup P_2$. If
$b_1a_3$ is red, then $a_2b_2b_1a_3xa_1y$ cover $V(G)$. Finally if
$b_1a_3$ is black, then $b_1a_3b_2x$ is a black $C_4$.

Therefore, $b_2y$ is black and $b_2x$ is blue. Consider $b_1a_3$. If
$b_1a_3$ is black, then $b_1a_3b_2y$ is a black $C_4$. If $b_1a_3$
is red, then $b_2b_1a_3xa_1y$ and $a_2$ cover more vertices than
$P_1\cup P_2$. If $b_1a_3$ is blue, then $b_1a_3a_2$ and $xa_1y$
cover more vertices than $P_1\cup P_2$.

\underline{Case 2.2:} $j>2$. Consider the edge $b_1b_j$. If $b_1b_j$
is blue, then $a_i\dots a_1 b_1 b_j$ plus a blue edge from $b_j$ to
$G'$ and $b_2 \dots b_{j-1}$ cover more vertices than $P_1\cup P_2$,
a contradiction. If $b_1b_j$ is red, then consider $b_2b_1b_j\dots
b_3$, a red path of length $j$. By Claim~\ref{ends}, there is a
cross-edge adjacent to two of $a_1,a_i,b_2,b_3$, and one of these
vertices, say $c$ (different from $b_1$) is non-adjacent to $G'$.
That is, $b_1xcy$ is a $C_4$ in $\overline G$, a contradiction. We
conclude $b_1b_j\in \overline{G}$. Now $a_ib_jb_1z$ is a path on $4$
vertices in $\overline G$, for any $z\in G'$. Therefore, any edge
$a_iz$, where $z\in G'$, is a red edge. If there is a red edge
$b_2z$, where $z\in G'$, then $b_1a_1\dots a_{i-1}$ and
$xa_izb_2\dots b_j$ cover more vertices than $P_1\cup P_2$, a
contradiction. Thus there is a blue edge $e$ from $b_2$ to $G'$. Now
consider the edge $b_2a_i$. If it were blue, then $b_1a_1\dots
a_ib_2$ extended with $e$ and $b_3\dots b_j$ would cover more
vertices than $P_1\cup P_2$, a contradiction. If $b_2a_i$ was red,
then $b_1a_1,\dots a_{i-1}$ and $xa_ib_2,\dots b_j$  would cover
more vertices than $P_1\cup P_2$, a contradiction. We conclude that
$b_2a_i\in \overline{G}$.

Next look at the pair $a_1,b_2$. It must be an edge $G$, otherwise
$a_1b_2a_ib_j$ is a $C_4$ in $\overline G$, a contradiction. If
$a_1b_2$ is red, then let $f$ be a red edge from $a_1$ to $U$, say
$f=a_1x$. Now $a_2\dots a_{i-1}$ and $ya_ixa_1b_2\dots b_j$ cover
more vertices than $P_1\cup P_2$, a contradiction. We conclude that
$a_1b_2$ is blue.

Consider the edge $a_ib_1$. If it is red, then $a_1\dots a_{i-1}$
and $xa_ib_1\dots b_j$ form a better pair. If $a_ib_1$ is blue, then
$b_1a_i\dots a_1b_2e$ and $b_3\dots b_j$ form a better pair. We
conclude $a_ib_1\in \overline{G}$. Now the $b_ja_ib_1z$ is a path on
4 vertices in $\overline{G}$, for any $z\in G'$. Therefore any
$b_jz$ in $\overline{G}$ would form a $C_4$. That is, all $b_jz$ are
blue edges.

Let $z$ be the endvertex of $e$ in $G'$. Now $a_i\dots a_1b_2zb_jx$
and $b_3\dots b_{j-1}$ cover more vertices than $P_1\cup P_2$,
giving a final contradiction. \qed\\

\noindent {\bf Acknowledgements.} We are indebted to an anonymous
referee for the extremely careful reading of the paper and many
helpful comments and suggestions.


\begin{thebibliography}{99}


\bibitem{A} P. Allen, Covering two-edge-coloured complete graphs with
two disjoint mon\-o\-chro\-mat\-ic cycles, {\em Combinatorics,
Probability and Computing,} {\bf 17} (2008), 471--486.

\bibitem{AY} J. Ayel, Sur l'existence de deux cycles supplementaires
unicolor\'es, disjoints et de couleurs differentes dans un graphe
complet bicolore, Thesis, University of Grenoble, (1979).

\bibitem{BBLS} F.S. Benevides, T. {\L}uczak, A. Scott, J. Skokan, M. White, Monochromatic cycles and the
monochromatic circumference in 2-coloured graphs,
{\em  Combin. Probab. Comput.} {\bf 21} (2012),  57–-87.



\bibitem{BT} S. Bessy, S. Thomass\'e, Partitioning a graph into a
cycle and an anticycle, a proof of Lehel's conjecture, {\em
Journal of Combinatorial Theory, Series B} {\bf 100}  (2010),
176--180.

\bibitem{EG} P. Erd\H os, T. Gallai, On maximal paths and circuits of
graphs, {\em Acta Math. Acad. Sci. Hungar.} {\bf 10} (1959),
337-356.

\bibitem{EGP} P. Erd\H os, A. Gy\'arf\'as, L. Pyber, Vertex
coverings by monochromatic cycles and trees, {\em Journal of
Combinatorial Theory, Series B} {\bf 51} (1991), 90--95.

\bibitem{FLPS} R. J. Faudree, S. L. Lawrence, T. D. Parsons, R. H. Schelp,
Path-Cycle Ramsey numbers, {\em Discrete Mathematics} {\bf 10}
(1974), 269--277.


\bibitem{GY} A. Gy\'arf\'as, Vertex coverings by monochromatic paths and cycles,
{\em Journal of Graph Theory} {\bf 7} (1983), 131--135.

\bibitem{rlogr} A. Gy\'arf\'as, M. Ruszink\'o, G. S\'ark\"ozy, E.
Szemer\'edi, An improved bound for the mono\-chromatic cycle
partition number, {\em Journal of Combinatorial Theory, Series B}
{\bf 96} (2006), 855--873.

\bibitem{GST} A. Gy\'arf\'as, G. Simonyi, \'A. T\'oth, Gallai
colorings and domination in multipartite digraphs,  {\em J. Graph Theory} {\bf 71} (2012),  278--292.


\bibitem{Gydisz} A. Gy\'arf\'as, Partition coverings and blocking sets
in hypergraphs (in Hungarian) {\em Communications of the Computer
and Automation Institute of the Hungarian Academy of Sciences} {\bf
71} (1977) 62 pp.

\bibitem{GYLSS} A. Gy\'arf\'as, J. Lehel, G. N. S\'ark\"ozy, R. H.
Schelp, Monochromatic Hamiltonian Berge cycles in colored complete
hypergraphs, {\em Journal of Combinatorial Theory, Ser. B} {\bf 98}
(2008), 342--358.

\bibitem{GRSSz} A. Gy\'arf\'as, M. Ruszink\'o, G. N. S\'ark\"ozy,
E. Szemer\'edi, Tripartite Ramsey numbers for paths, {\em Journal of
Graph Theory} {\bf 55} (2007), 164--170.

\bibitem{GYRSSZ}
A. Gy\'arf\'as, M. Ruszink\'o, G. S\'ark\"ozy, E. Szemer\'edi,
Three-color Ramsey numbers for paths, {\em Combinatorica} {\bf 27}
(2007), 35--69.

\bibitem{GRSSZ}
A. Gy\'arf\'as, M. Ruszink\'o, G. S\'ark\"ozy, E. Szemer\'edi,
Partitioning 3-colored complete graphs into three monochromatic
cycles, {\em Electronic Journal of Combinatorics} {\bf 18} (2011),
Nr. 53.

\bibitem{GYS1}
A. Gy\'arf\'as, G. N. S\'ark\"ozy, Gallai colorings of non-complete
graphs, {\em Discrete Mathematics} {\bf 310} (2010), 977--980.

\bibitem{GYS2}
A. Gy\'arf\'as, G. N. S\'ark\"ozy, Star versus two stripes Ramsey
numbers and a conjecture of Schelp,  {\em Combin. Probab. Comput.} {\bf 21} (2012),  179--186.
\bibitem{KS}
J. Koml\'os and M. Simonovits, Szemer\'edi's Regularity Lemma and
its applications in Graph Theory, in Combinatorics, Paul Erd\H{o}s
is Eighty (D. Mikl\'{o}s, V.T. S\'os and T. Sz\H onyi, Eds.), pp.
295-352, Bolyai Society Mathematical Studies, Vol. 2, Budapest,
1996.

\bibitem{KOV}
T. K\H ov\'ari, V. T. S\'os, P. Tur\'an, On a problem of K.
Zarankiewicz, {\em Colloq. Math.} {\bf 3} (1954) 50--57.


\bibitem{LO}
L. Lov\'asz, Combinatorial problems and exercises, 2nd Edition,
American Mathematical Society, Providence, Rhode Island, 2007.


\bibitem{L}
T. \L uczak, $R(C_n,C_n,C_n)\le(4+o(1))n$, {\em Journal of
Combinatorial Theory, Ser. B} {\bf 75} (1999), 174--187.

\bibitem{LRS}
T. \L uczak, V. R\"{o}dl, E. Szemer\'edi, Partitioning two-colored
complete graphs into two monochromatic cycles, {\em Probability,
Combinatorics and Computing} {\bf 7} (1998), 423--436.

\bibitem{posa}
L. P\'osa, On the circuits of finite graphs, {\em MTA Mat. Kut. Int.
K\"ozl.} {\bf 8} (1963), 355--361.

\bibitem{PR}
A. Pokrovskiy, Partitioning $3$-coloured complete graphs into three
monochromatic paths, {\em Electronic Notes in Discrete Mathematics}
{\bf 38} (2011) 717--722.

\bibitem{PR1}
A. Pokrovskiy, Partitioning edge-coloured complete graphs into
monochromatic cycles and paths, {\em arXiv:1205.5492v1}

\bibitem{SG}
G. N. S\'ark\"ozy, On $2$-factors with $k$ components, {\em Discrete
Mathematics} {\bf 308} (2008), 1962--1972.

\bibitem{SA}
G. N. S\'ark\"ozy, Monochromatic cycle partitions of edge-colored
graphs, {\em Journal of graph theory} {\bf 66} (2011), 57--64.

\bibitem{SCH}
R. H. Schelp, A minimum degree condition on a Ramsey graph which
arrows a path, submitted to Discrete Mathematics.

\bibitem{Sz}
E. Szemer\'edi, Regular partitions of graphs, Colloques
Internationaux C.N.R.S. $\mbox{N}^{\underline{o}}$ 260 - Probl\`emes
Combinatoires et Th\'eorie des Graphes, Orsay (1976), 399--401.

\end{thebibliography}
\end{document}